\documentclass[12pt]{amsart}

\newcommand{\lrw}{\longrightarrow}
\newcommand{\Map}{\longmapsto}
\newcommand{\beqa}{\begin{eqnarray*}}
\newcommand{\eeqa}{\end{eqnarray*}}

\renewcommand{\a}{\mbox{$\alpha$}}

\newcommand{\hs}{\hspace{.2in}}

\newcommand{\Ad}{{\rm Ad}}

\def \fg{{\mathfrak g}}
\def \fh{{\mathfrak h}}
\def \fk{{\mathfrak k}}
\def \fp{{\mathfrak p}}
\def \fq{{\mathfrak q}}
\def \fa{{\mathfrak a}}
\def \fb{{\mathfrak b}}
\def \fn{{\mathfrak n}}
\def \ft{{\mathfrak t}}

\def \sl{{\mathfrak s}{\mathfrak l}}

\newcommand{\fu}{\mathfrak u}

\def \R{{\mathbb R}}
\def \P{{\mathbb P}}
\def \C{{\mathbb C}}

\def \i{{i}}
\def \fd{{\bullet}}
\def \proof{{\noindent {\bf Proof.} \ \ }}
\renewcommand{\qed}{\begin{flushright} {\bf Q.E.D.}\ \ \ \ \
                  \end{flushright} }

\newtheorem{thm}{Theorem}[section]
\newtheorem{lem}[thm]{Lemma}
\newtheorem{prop}[thm]{Proposition}
\newtheorem{proposition}[thm]{Proposition}
\newtheorem{cor}[thm]{Corollary}
\newtheorem{rem}[thm]{Remark}
\newtheorem{rems}[thm]{Remarks}
\newtheorem{exam}[thm]{Example}

\newtheorem{dfn-nota}[thm]{Definition-Notation}
\newtheorem{dfn-prop}[thm]{Definition-Proposition}

\numberwithin{equation}{section}

\begin{document}

\baselineskip=16pt

\title[A Poisson structure on symmetric spaces]{A Poisson structure on  
compact symmetric spaces}

\author[P. Foth and J.-H. Lu]{Philip Foth and Jiang-Hua Lu}

\address{Department of Mathematics, University of Arizona, Tucson, AZ 85721-0089
\newline \indent Department of Mathematics, University of Hong Kong, Pokfulam Road, Hong Kong}

\email{foth@math.arizona.edu, jhlu@maths.hku.hk}



\date{September 20, 2003}

\begin{abstract}
We present some basic results on a natural Poisson structure on any compact 
symmetric space.  
The symplectic leaves of this structure are related to the orbits of the
corresponding real semisimple group on the complex flag manifold. 
\end{abstract}

\maketitle

\section{Introduction and the Poisson structure $\pi_0$ on $U/K_0$.}

Let $\fg_0$ be a real semi-simple Lie algebra, and let
$\fg$ be its complexification. Fix a Cartan 
decomposition $\fg_0 = \fk_0 + \fp_0$ of $\fg_0$, and let $\fu$
be the compact real form of $\fg$ given by $\fu = \fk_0 + i\fp_0$.
Let $G$ be the connected and simply
connected Lie group with Lie algebra $\fg$, and let $G_0, K_0$, and $U$ be the
connected subgroups of $G$ with Lie algebras $\fg_0, \fk_0$, and $\fu$ 
respectively. Then $K_0 = G_0 \cap U$, 
and $U/K_0$ is the compact dual of the non-compact
Riemannian symmetric space $G_0/K_0$. In this paper, we will define
a Poisson structure $\pi_0$ on $U/K_0$ and study some of its properties.

The definition of $\pi_0$ depends on a choice of an {\it Iwasawa--Borel} subalgebra
of $\fg$ relative to $\fg_0$. Recall \cite{huckle-wolf} that a Borel subalgebra $\fb$
of $\fg$
is said to be Iwasawa relative to $\fg_0$ 
if $\fb \supset \fa_0 + \fn_0$ for some Iwasawa decomposition $\fg_0 =
\fk_0 + \fa_0 + \fn_0$ of $\fg_0$. Let $Y$ be the variety of all Borel subalgebras of 
$\fg$. Then $G$ acts transitively
on $Y$ by conjugations, and   $\fb \in Y$ is  
Iwasawa relative to $\fg_0$ if and only if it lies in the unique closed orbit of
$G_0$ on $Y$ \cite{huckle-wolf}. Denote by
$\tau$ and $\theta$ the complex conjugations on $\fg$ 
with respect to $\fg_0$ and $\fu$ 
respectively. 
Throughout this paper,
we will fix an Iwasawa--Borel subalgebra $\fb$ relative to $\fg_0$ and 
a Cartan subalgebra $\fh \subset
\fb$ of $\fg$ that is stable under both $\tau$ and $\theta$. Let $\Delta^+$ be the 
set of roots for $\fh$ determined by $\fb$, and let $\fn$ be the complex
span of root vectors for roots in $\Delta^+$, so that $\fb = \fh + \fn$.
Let $\fa = \{x \in \fh: \theta(x) = -x\}$.
Let $\fa_0 = \fa \cap \fg_0$ and 
$\fn_0 = \fn \cap \fg_0$. Then $\fg_0 = \fk_0 + \fa_0 + \fn_0$
is an Iwasawa decomposition of $\fg_0$. 

We can define a Poisson structure $\pi_0$ on $U/K_0$ as follows:
let $\ll \, , \, \gg$ be the Killing form of
$\fg$. For each $\alpha \in \Delta^+$, choose a root vector $E_\alpha$ such that
$\ll E_\alpha, \, \theta(E_{\alpha}) \gg = -1$. Let $E_{-\alpha} = 
-\theta(E_\alpha)$, and let
$X_\alpha = E_\alpha - E_{-\alpha}$ and $Y_\alpha =i(E_\alpha + E_{-\alpha})$. 
Then
$X_\alpha, Y_\alpha \in \fu$ for each $\alpha \in \Delta^+$. Set
$$
\Lambda=\frac{1}{4}\sum_{\alpha \in \Delta^+}X_{\alpha}\wedge 
Y_{\alpha}\in\fu\wedge\fu,
$$
and define the bi-vector field $\pi_U$ on $U$ by
$$
\pi_U=\Lambda^r-\Lambda^l,
$$ where  
 $\Lambda^r$ and $\Lambda^l$ are respectively the right and left invariant  
bi-vector fields on $U$
with value $\Lambda$ at the identity element. Then $\pi_U$ is a 
Poisson bivector field, and  
$(U,\pi_U)$ is the Poisson-Lie group defined by
the Manin triple $(\fg, \fu, \fa + \fn)$ \cite{LW}. 

The group $G$ acts on $U$ from the right via $u^g=u_1$, if $ug=bu_1$  
for some $b\in AN$, where $A = \exp \fa$ and $N = \exp \fn$.  
Therefore every subgroup of $G$, for example $AN$ or $G_0$, also acts on $U$. The   
symplectic leaves of $\pi_U$ are precisely the orbits of the right $AN$-action.  
These leaves are parameterized by the torus $T = \exp(i\fa)$ and the Weyl group $W$
of $(U, \fh)$.   
The Poisson structure $\pi_U$ is both left and right $T$-invariant, and it   
descends to the so-called
Bruhat Poisson structure on $T\backslash U$, whose symplectic leaves are
precisely the Bruhat cells of $T \backslash U \cong B\backslash G$ as
the orbits of the Borel group $B=TAN$. We refer to \cite{LW}
for details. 

\begin{prop}
\label{prop_def-piz}
There exists a Poisson structure $\pi_0$ on $U/K_0$  
such that the natural projection
$p: (U, \pi_U) \to (U/K_0, \pi_0)$ is a Poisson map. The symplectic leaves of 
the Poisson  
structure $\pi_0$ are precisely the  
projections of the $G_0$-orbits on $U$ via the map $p$.     
\end{prop}      

\proof To show that the Poisson structure $\pi_U$ descends to the  
quotient $U/K_0$, it is enough to show that  
the annihilator space $\fk_0^{\perp}$  of $\fk_0$ inside $\fu^*$, which
is identified with $\fa + \fn$, is a Lie subalgebra of  
$\fa+\fn$. The bilinear form which is used in this identification is  
the imaginary part of the Killing
form $\ll \, , \, \gg$ of $\fg$. We observe that being a real form of $\fg$,
$\fg_0$ is isotropic with respect to ${\rm Im} \ll \, , \, \gg$, 
which implies that $\fk_{0}^{\perp} \subset \fa_0 + \fn_0$. It then follows
for dimension reason that $\fk_{0}^{\perp} = \fa_0 + \fn_0$, which is 
a Lie subalgebra of $\fa + \fn$.

For the statement concerning the symplectic leaves of $\pi_0$,  
we observe that $(X, \pi_0)$ is a $(U, \pi_U)$-Poisson  
homogeneous space, and then apply  
\cite[Theorem 7.2]{LuDuke}. \qed

\begin{rem} {\em For the case when the Satake diagram of $\fg_0$ has no black 
dots, the Poisson structure $\pi_0$
was considered by Fernandes in \cite{Fernandes}. }
\end{rem}

In this paper, we will study some properties of 
the symplectic leaves
of $\pi_0$. Recall that $Y$ is the variety of all Borel subalgebras of $\fg$. 
We will show that the set of symplectic leaves of $\pi_0$ is essentially
parameterized by the set of $G_0$-orbits in $Y$,
which have been studied extensively because of their importance in
the representation theory of $G_0$. More precisely, 
let
$q: U \to Y$ be surjective map $u \mapsto \Ad_{u}^{-1} \fb \in Y$.
Then the map ${\mathcal O} \mapsto
p(q^{-1}({\mathcal O}))$ gives a bijective correspondence between the set of
$G_0$-orbits in $Y$ and the set of $T$-orbits of symplectic leaves in
$U/K_0$. In particular, there are finitely many families of symplectic leaves.
In each family leaves are translates of one another by elements in $T$. 
Moreover, $\pi_0$ has open symplectic leaves
if and only if $\fg_0$ has a compact Cartan subalgebra, in which case, 
the number of open symplectic leaves is the same as the number of
open $G_0$-orbits in $Y$, and 
each open symplectic leaf is diffeomorphic to $G_0/K_0$.
When $X$ is Hermitian symmetric, the Poisson  
structure $\pi_0$ is shown to be 
the sum of the Bruhat Poisson structure \cite{LW} and
a multiple of any non-degenerate invariant Poisson structure.   

We also show that the $U$-invariant Poisson 
cohomology $H^{\fd}_{\pi_0, U}(U/K_0)$ is isomorphic to the De Rham 
cohomology of $U/K_0$. 
The full Poisson cohomology and some further properties of $\pi_0$ will be
studied in a future paper. 

Throughout the paper, if $Z$ is a set and if $\sigma$ is an 
involution on $Z$, we will use $Z^\sigma$ to denote the fixed point
set of $\sigma$ in $Z$.

\section{Symplectic leaves of $\pi_0$ and $G_0$-orbits in $Y$.}
\label{sec_5}

By Proposition \ref{prop_def-piz}, symplectic leaves of $\pi_0$ 
are precisely  
the projections to $U/K_0$ of $G_0$-orbits
in $U$. Here, recall that $G_0$ acts on $U$
as a subgroup of $G$, and $G$ acts on $U$ from the right by
\begin{equation}
\label{eq_G-on-U}
u^g = u_1, \hspace{.2in} {\rm if} \hspace{.2in} ug = b u_1 \hspace{.1in} 
{\rm for}
\hspace{.1in} b \in AN,
\end{equation}
where $u \in U$ and $g \in G$.  It is easy to see
that the above right action
of $G$ on $U$ descends to an action of $G$ on $T \backslash U$. 
On the other hand,
the map $U \to Y: u \mapsto \Ad_{u}^{-1} \fb$ gives a $G$-equivariant
identification of $Y$ with $T \backslash U$. This identification
will be used throughout the paper.
The $G_0$-orbits on $Y$ have been studied extensively (see,
for example, \cite{Matsuki:cartan} and \cite{wolf1}). In particular,  
there are finitely many $G_0$-orbits in $Y$. We will now formulate
a precise connection between symplectic leaves of $\pi_0$ and $G_0$-orbits 
in $Y$.

Let $X = U/K_0$. For $x \in X$, let $L_x$ be the symplectic leaf of $\pi_0$  
through $x$. Since $T$ acts by Poisson diffeomorphisms, for each
$t \in T$, the set $tL_x = \{tx_1: x_1 \in L_x\}$ is again a  
symplectic leaf of $\pi_0$. Let
\[
{\mathcal S}_x = \bigcup_{t \in T} t L_x \subset X.
\]
For $y \in Y$, let ${\mathcal O}_y$ be the $G_0$-orbit
in $Y$ through $y$. Let $p: U \to X= U/K_0$ and
$q: U \to  Y =T \backslash U $ be the natural projections.

\begin{prop}
\label{prop_cor1}
Let $x \in X$ and $y \in Y$ be such that $p^{-1}(x) \cap q^{-1}(y)
\neq \emptyset$. Then
\[
p(q^{-1}({\mathcal O}_y)) = {\mathcal S}_x, \hs {\rm and}
\hs q(p^{-1}({\mathcal S}_x)) = {\mathcal O}_y.
\]
\end{prop}

\proof Let $u \in p^{-1}(x) \cap q^{-1}(y)$, and let
$u^{G_0}$ be the $G_0$-orbit in $U$ through $u$. It is  
easy to show that  
\[
q^{-1}({\mathcal O}_y) = p^{-1}({\mathcal S}_x) = \bigcup_{t \in T}
t(u^{G_0}).
\]
Thus,
\[
p(q^{-1}({\mathcal O}_y)) = \bigcup_{t \in T} t p(u^{G_0}) = {\mathcal S}_x,
\]
and
\[
q(p^{-1}({\mathcal S}_x)) = q(u^{G_0}) = {\mathcal O}_y.
\]
\qed

\begin{cor}
\label{cor_bijection}
Let ${\mathcal O}_Y$ be the collection of $G_0$-orbits in $Y$, and
let ${\mathcal S}_X$ be the collection of all the subsets
${\mathcal S}_x, x \in X$. Then the map
\[
{\mathcal O}_Y \longrightarrow {\mathcal S}_X: \ \   
{\mathcal O} \longmapsto
p(q^{-1}({\mathcal O}))
\]
is a bijection with the inverse given by ${\mathcal S} \mapsto  
q(p^{-1}({\mathcal S}))$.
\end{cor}

We now recall some facts about $G_0$-orbits in $Y$ 
from \cite{R-S:orbits} which we will use to
compute the dimensions of symplectic leaves of $\pi_0$. 
Since  \cite{R-S:orbits} is based on the choice of a
Borel subalgebra in an open $G_0$-orbit in $Y$, we will
restate the relevant results from \cite{R-S:orbits} in
Proposition \ref{prop_RS-results} to fit our set-up.

Let $\ft = i \fa$ be the Lie algebra of $T$, and let $N_U(\ft)$ be the
normalizer subgroup of $\ft$ in $U$. Set
\[
{{\mathcal V}} = \{u \in U: \, u \tau(u)^{-1} \in N_U(\ft)\}.
\]
Then $u \in {{\mathcal V}}$ if and only if 
$\Ad_{u}^{-1} \fh$ is $\tau$-stable. Clearly 
${{\mathcal V}}$ is invariant under the left translations
by elements in $T$ and the right translations by elements in $K_0$. Set
\[
V = T\backslash {{\mathcal V}}/K_0.
\]
Then we have a well-defined map
\[
 V \lrw {{\mathcal O}}_Y: \, \, v \Map {{\mathcal O}}(v),
\]
 where for $v = TuK_0 \in V$, ${{\mathcal O}}(v)$ is the $G_0$-orbit in $Y$ 
through the point
$\Ad_{u}^{-1} \fb \in Y$. Let $W = N_U(\ft)/T$ be the Weyl group. 
Then we also have 
the well-defined map
\[
\psi: \, \, V \lrw W: \, \, v = TuK_0 \Map u\tau(u)^{-1}T \in W.
\]
For $w \in W$, let $l(w)$ be the length of $w$. 

\begin{prop}
\label{prop_RS-results}
1) The map $v \mapsto {{\mathcal O}}(v)$ is a bijection between the set $V$ and the 
set ${{\mathcal O}}_Y$ of all $G_0$-orbits in $Y$;

2) For $v \in V$, the co-dimension of ${{\mathcal O}}(v)$ in $Y$ is equal to
$l(\psi(v)w_bw_0)$, where $w_0$ is the longest 
element of $W$, and $w_b$ is the longest element of the
subgroup of $W$ generated by the black dots of the Satake diagram of $\fg_0$.
\end{prop}

\begin{rems} 
\label{rem_satake-diagram}
{\em 1) Since $\tau$ leaves $\fa$ invariant, it acts on the set of roots for $\fh$
by $(\tau \alpha)(x) = \alpha(\tau(x))$ for $x \in \fa$. We know from \cite{araki} 
that the black dots in the Satake diagram of $\fg_0$ correspond precisely to the
simple roots $\alpha$ in $\Delta^+$ such that $\tau(\alpha) = -\alpha$. 
Moreover, if $\alpha \in \Delta^+$ and if $\tau (\alpha) \neq -\alpha$, then
$\tau(\alpha) \in \Delta^+$; 

2) We now point out how Proposition \ref{prop_RS-results} follows from
results in \cite{R-S:orbits}. 
 Let $u_0 \in U$ be such that
$\fb^\prime :=\Ad_{u_0} \fb$ lies in an open $G_0$-orbit in $Y$ and 
$\fh^\prime: =\Ad_{u_{0}} \fh$ is $\tau$-stable. The pair $(\fg_0, \fb^\prime)$ 
is called
a {\it standard pair} in the terminology of \cite[No.1.2]{R-S:orbits}.  Let
$\ft^\prime = \Ad_{u_0} \ft$, $T^\prime = u_0Tu_{0}^{-1}$, 
and  $N_U(\ft^\prime) = u_{0} N_U(\ft) u_{0}^{-1}$.
Let
\[
{{\mathcal V}}^\prime = \{u^\prime \in U: \, \, u^\prime
 \tau(u^\prime)^{-1} \in N_U(\ft^\prime)\},
\]
and let $V^\prime = T^\prime \backslash {{\mathcal V}}^\prime /K_0$. For
$v^\prime = T^\prime u^\prime K_0$, let
${{\mathcal O}}(v^\prime)$ be the $G_0$-orbit in $Y$ through the point
$\Ad_{u^\prime}^{-1} \fb^\prime \in Y$. 
Then 
\cite[Theorem 6.1.4(3)]{R-S:orbits} says that
the map $V^\prime \to {{\mathcal O}}_Y: v^\prime \to {{\mathcal O}}(v^\prime)$
is a bijection between the set $V^\prime$ and the set ${{\mathcal O}}_Y$ of 
$G_0$-orbits in $Y$, and 
\cite[Theorem 6.4.2]{R-S:orbits} says that 
the co-dimension of ${{\mathcal O}}(v^\prime)$ in $Y$ is
the length of the element $\phi(v^\prime)$ 
in the Weyl group $W^\prime =N_U(\ft^\prime)/T^\prime $
defined by $u^\prime \tau(u^\prime)^{-1} \in N_U(\ft^\prime)$.
Since $\fb = \Ad_{u_0}^{-1} \fb^\prime$ lies in the 
unique closed $G_0$-orbit
in $Y$, it follows from \cite[No. 1.6]{R-S:orbits} that 
$u_0 \tau(u_0)^{-1} \in N_U(\ft^\prime)$ defines the element in $W^\prime$ 
that corresponds to $w_bw_0 \in W$ under the natural
identification of $W$ and $W^\prime$. It is also easy to see that
${{\mathcal V}}^\prime = u_0 {{\mathcal V}}$, and if 
$v^\prime = T^\prime u^\prime
K_0 \in V^\prime$ and $v = T(u_{0}^{-1} u^\prime)K_0 \in V$
for
$u^\prime \in {{\mathcal V}}^\prime$,  then
${{\mathcal O}}(v^\prime) = {{\mathcal O}}(v)$, and
$\phi(v^\prime) \in W^\prime$ corresponds to  $\psi(v)w_bw_0 \in W$
under the natural
identification of $W$ and $W^\prime$. 
It is now clear that Proposition \ref{prop_RS-results} holds.
Statement 2) of Proposition \ref{prop_RS-results} can also be seen
directly from Lemma \ref{lem_diffeo} below;

3) Starting from a complete collection of 
representatives of equivalence
classes of strongly orthogonal real roots for the Cartan subalgebra
$\fh^\tau$ of $\fg_0$, it is possible, by using Cayley transforms,
to explicitly construct a set of representatives of $V$ in 
${{\mathcal V}}$. This is done in \cite[Theorem 3]{Matsuki:cartan}.

4) The three involutions $\tau, w_0$ and $w_b$ on 
$\Delta = \Delta^+ \cup (-\Delta^+)$ commute with each other. Indeed, 
since $\tau$ commutes with the reflection defined by every black
dot on the Satake diagram, $\tau$ commutes with $w_b$.
We know from Remark (\ref{rem_satake-diagram}) that $\tau w_b(\Delta^+) = 
\Delta^+$, so $\tau w_b$ defines an automorphism
of the Dynkin diagram of $\fg$. It is well-known that $-w_0$ is
in the center of the group of all automorphisms of the Dynkin diagram of
$\fg$ (this can be checked, for example, case by case). Thus $w_0$ commutes with 
$\tau w_b$. To see that $w_0$ commutes with $w_b$, note by directly
checking case by case that $-w_0$ maps a simple black root 
on the Satake diagram of $\fg_0$ to another such  simple black root. 
Thus $w_0w_bw_0$
is still in the subgroup $W_b$ of $W$ generated by
the set of all black simple roots. It follows that 
$w_0w_b$ and $w_b w_0 = w_0(w_0w_bw_0)$ are in the same right $W_b$ coset in $W$.
Since $l(w_0w_b) = l(w_bw_0) = l(w_0)-l(w_b)$, we know that $w_0w_b=w_bw_0$
by the uniqueness of minimal length representatives of right $W_b$ cosets
in $W$. Thus $w_0$ commutes with both $\tau$ and $w_b$. 
These remarks will be used
in the proof of Lemma \ref{lem_diffeo}.
}
\end{rems}

\section{Symplectic leaves of $\pi_0$.}

Recall that $p: U \to U/K_0$ and $q: U \to Y = 
T\backslash U$
are the natural projections. For each
$v \in V = T\backslash {{\mathcal V}}/K_0$, set 
\[
{{\mathcal S}}(v) = p(q^{-1}({{\mathcal O}}(v))) \subset 
U/K_0.
\]
By Corollary \ref{cor_bijection}, we have a disjoint union
\[
U/K_0 = \bigcup_{v \in V}  {{\mathcal S}}(v).
\]
Moreover, each ${{\mathcal S}}(v)$ is a union of symplectic leaves
of $\pi_0$, all of which are 
 translates of each other by elements in $T$. Thus it is 
enough to understand  one single leaf in ${{\mathcal S}}(v)$.
Recall that
$G$ acts on $U$ from the right by $(u, g) \mapsto u^g$ as described 
in (\ref{eq_G-on-U}). 

\begin{lem}
\label{lem_leaves-double-cosets}
For every $u \in U$, the map
\[
(G_0 \cap u^{-1}(AN)u) \backslash G_0/K_0 \lrw U/K_0:\, \, 
(G_0 \cap u^{-1}(AN)u) g_0K_0 \Map u^{g_0}K_0, \hs g_0 \in G_0,
\]
gives a diffeomorphism between the double coset space 
$(G_0 \cap u^{-1}(AN)u) \backslash G_0/K_0$ and the symplectic leaf of $\pi_0$ 
through the point $uK_0 \in U/K_0$.
\end{lem}

\proof 
Consider the $G_0$-action on $U$ as a subgroup of $G$. By (\ref{eq_G-on-U}), 
the induced action of $K_0$ on $U$ is by left translations.
It is easy to see that the stabilizer subgroup of $G_0$ at $u$ is
$G_0 \cap u^{-1}(AN)u$. 
Let  $u^{G_0}$ be the $G_0$-orbit in $U$ 
through $u$. Then 
\[
u^{G_0} \cong (G_0 \cap u^{-1}(AN)u) \backslash G_0.
\]
Since the action of $K_0$ on $u^{G_0}$ by left translations is free, we see that the
double coset space $(G_0 \cap u^{-1}(AN)u) \backslash G_0/K_0$ is smooth. 
Lemma \ref{lem_leaves-double-cosets} now follows from 
Proposition \ref{prop_def-piz}.
\qed
 
Assume now that $u \in {{\mathcal V}}$. To better understand
the group
$G_0 \cap u^{-1}(AN)u$, we introduce the involution $\tau_u$ on $\fg$:
\[
\tau_u = \Ad_u \tau \Ad_{u}^{-1}= 
\Ad_{u\tau(u^{-1})} \tau: \, \, \fg \lrw \fg.
\]
The fixed point set of $\tau_u$ in $\fg$ is the real form 
$\Ad_u \fg_0$ of $\fg$. We will use the same letter for the lifting 
of $\tau_u$ to $G$. Since $\tau_u$
leaves $\fa$ invariant, it acts on the set of roots for $\fh$ by
$(\tau_u \alpha)(x) = \alpha(\tau_u(x))$ for $x \in \fa$. Recall that
associated to $v=TuK_0 \in V$ we have the Weyl group element $\psi(v)w_bw_0$.
Let 
\[
N_v = 
N \cap (\dot{w}N^- \dot{w}^{-1}), 
\]
where $\dot{w} \in U$
is any representative of $\psi(v)w_bw_0 \in W$.

\begin{lem}
\label{lem_diffeo}
For  any $u \in {{\mathcal V}}$ and $v = T u K_0 \in V$, 

1)  $\Delta^+ \cap \tau_u(\Delta^+) = \Delta^+ \cap (\psi(v)w_bw_0)
(-\Delta^+)$;

2) $N_v$ is $\tau_u$-invariant and $G_0 \cap u^{-1} N u 
 =  u^{-1}(N_{v})^{\tau_u}u 
=(u^{-1} N_v u)^\tau$
is connected;

3) the map
\begin{equation}
\label{eq_3}
M: \, \, (G_0 \cap u^{-1}Tu) \times (G_0 \cap u^{-1}Au) \times
(G_0 \cap u^{-1}Nu) \lrw G_0 \cap u^{-1}(TAN)u
\end{equation}
given by $M(g_1, g_2, g_3)  =g_1g_2g_3$ 
is a diffeomorphism. 
\end{lem}

\proof 
1) Recall that $\psi(v) \in W$ is the element defined by $u \tau(u)^{-1} \in N_U(\ft)$.
Then $\tau_u(\alpha) = \psi(v) \tau(\alpha)$ for every $\alpha \in \Delta$.
Thus $\tau_u(\alpha) \in \Delta^+$ if and only if $\psi(v)\tau(\alpha) \in \Delta^+$,
which is in turn equivalent to $w_0\tau w_b \psi(v) \tau(\alpha) \in - \Delta^+$
because $w_0 \tau w_b(\Delta^+) = -\Delta^+$. Since the three involutions
$w_0, \tau$ and $w_b$ commute with each other by Remark \ref{rem_satake-diagram}, 
we have $w_0\tau w_b \psi(v) \tau = (\psi(v) w_bw_0)^{-1}$. This proves 1).

2) We know from 1) that $\Delta^+ \cap (\psi(v)w_bw_0)
(-\Delta^+)$ is $\tau_u$-invariant. Thus $N_v$ is $\tau_u$-invariant. 
Clearly $u^{-1}(N_{v})^{\tau_u}u \subset G_0 \cap u^{-1} N u$. 
Let $N_{v}^{\prime} = N \cap \dot{w}N\dot{w}^{-1}$. Then
$N = N_vN_{v}^{\prime}$ is a direct product, and we know from
1) that $\tau_u(N_{v}^{\prime}) \subset N^{-}$. Suppose now that $n \in N$
is such that 
$u^{-1}n u \in G_0 \cap u^{-1} N u$.
Write $n = m m^\prime$ with $m \in N_v$ and $m^\prime \in N_{v}^{\prime}$.
Then from $\tau_u(n) = n$ we get
$\tau_u(m^\prime) = \tau_u(m^{-1}) n \in N^{-} \cap N 
= \{e\}$. Thus $m^\prime = e$, and  $n = m \in 
(N_v)^{\tau_u}$. Since the exponential map for the group $u^{-1}(AN)u$ is a diffeomorphism,
$(u^{-1}(AN)u)^\tau$ is the connected subgroup of $u^{-1}(AN)u$ with
Lie algebra $(\Ad_{u}^{-1}(\fa + \fn))^\tau$. This shows 2).

We now prove 3).
Since $\Ad_{u}^{-1} \fh$ is $\tau$-invariant, 
the Lie algebra $\fg_0 \cap \Ad_{u}^{-1}\fb$ of $G_0 \cap u^{-1}(TAN)u$
is the direct sum of
the Lie algebras of the three subgroups on the left hand side of (\ref{eq_3}).
Thus the map $M$ is a local diffeomorphism. It is also easy to see that
$M$ is one-to-one. Thus it remains to show that $M$ is onto. 
Suppose that 
$h \in TA$ and $n \in N$ are such that $u^{-1} (hn) u \in G_0$. Then
$\tau_u(hn) = hn.$ 
Write $n = m m^\prime$ with $m \in N_v$ and $m^\prime \in N_{v}^{\prime}$.
Then from $\tau_u(hn) = hn$ we get
$\tau_u(m^\prime) = \tau_u(m^{-1}) \tau_u(h^{-1}) hn \in N^{-} \cap
HN = \{e\}$. Thus $m^\prime = e$, and $\tau_u(h) = h$ and $n = m \in 
(N_v)^{\tau_u}$. 
 If $h = ta$ with $t \in T$ 
and $a \in A$, 
it is also easy to see that 
$\tau(h) = h$ implies that $\tau_u(t) = t$ and $\tau_u(a) = a$. 
 \qed
 
In particular, we see that  $G_0\cap u^{-1}(AN)u$ is a contractible subgroup 
of $G_0$. Since Lemma \ref{lem_leaves-double-cosets} states that
the symplectic leaf of $\pi_0$ through the point $uK_0$ is diffeomorphic to 
$(G_0 \cap u^{-1}(AN)u) \backslash G_0/K_0$, we see that this leaf is the base
space of a smooth fibration with contractible total space and fiber. 
Thus we have:   
\begin{proposition} 
Each symplectic leaf of the Poisson structure $\pi_0$ is contractible.
\end{proposition}

\begin{rem}
\label{rem_another}
{\rm
Since $\dim(Y) = \dim ((G_0 \cap u^{-1}(TA) u)\backslash G_0)$, it is also
clear from 3) of Lemma \ref{lem_diffeo} that the codimension of 
${{\mathcal O}}(v)$ in $Y$ is $l(\psi(v)w_bw_0)$. See Proposition 
\ref{prop_RS-results}.
}
\end{rem}

It is a basic 
fact \cite{wolf1} that associated to each $G_0$-orbit in $Y$ there is a 
unique  $G_0$-conjugacy class of $\tau$-stable Cartan subalgebras of $\fg$.
For $u \in {{\mathcal V}}$ and $v = TuK_0 \in V$, the 
$G_0$-conjugacy class of $\tau$-stable Cartan subalgebras of $\fg$
associated to  ${{\mathcal O}}(v)$ is that defined by 
$\Ad_{u}^{-1} \fh$. The intersection $(\Ad_{u}^{-1} \fh) \cap \fg_0$
is a Cartan subalgebra of $\fg_0$.  
Regard both $\tau$ and $\psi(v)$ as maps 
on $\fh$ so that $\psi(v)\tau = \tau_u|_{{\mathfrak h}}:\fh \to \fh$. 
Then we have    
\[
(\Ad_{u}^{-1} \fh) \cap \fg_0
=(\Ad_{u}^{-1} \fh)^\tau = \Ad_{u}^{-1}(\fh^{\psi(v)\tau}).
\]
Since $\psi(v) \tau$ commutes with $\theta$, 
it leaves both $\ft = \fh^{\theta}$ and $\fa = \fh^{-\theta}$ 
invariant, and we have
\[
(\Ad_{u}^{-1} \fh) \cap \fg_0 = \Ad_{u}^{-1}(\ft^{\psi(v)\tau} + \fa^{\psi(v)\tau}).
\]
The subspaces $\Ad_{u}^{-1}(\ft^{\psi(v)\tau})$ and
$\Ad_{u}^{-1}( \fa^{\psi(v)\tau})$ are respectively the toral and vector parts of the 
Cartan subalgebra $(\Ad_{u}^{-1} \fh) \cap \fg_0$ of $\fg_0$. Set
\begin{eqnarray}
\label{eq_tv}
t(v) &=& \dim (\ft^{\psi(v)\tau}) =\dim (\Ad_{u}^{-1}(\ft^{\psi(v)\tau}))
=\dim(G_0 \cap u^{-1}T u)\\
\label{eq_av}
a(v) &= &\dim (\ft^{\psi(v)\tau}) =\dim (\Ad_{u}^{-1}(\fa^{\psi(v)\tau}))=
\dim (G_0 \cap u^{-1}A u).
\end{eqnarray}

\begin{thm}
\label{thm:83} 
For every $v \in V$,

1) every symplectic leaf $L$ in ${{\mathcal S}}(v)$  has  
dimension 
\[
\dim L = \dim ({{\mathcal O}}(v)) -\dim (K_0) + t(v),
\]
so the co-dimension of $L$ in $U/K_0$ is $a(v) + l(\psi(v)w_bw_0)$;

2) the family of
symplectic leaves in ${{\mathcal S}}(v)$ is  parameterized by the quotient
torus  $T/T^{\psi(v)\tau}$.
\end{thm}   

\proof Let $u$ be a representative of $v$ in ${{\mathcal V}} \subset U$.
Let $x = uK_0 \in U/K_0$, and let
$L_x$ be the symplectic leaf of $\pi_0$ through $x$. We only need to compute
the dimension of $L_x$.
 Let  
$u^{G_0}$ be the $G_0$-orbit in $U$ 
through $u$. We know from Lemma \ref{lem_leaves-double-cosets}
that $u^{G_0} \cong (G_0 \cap u^{-1}(AN)u) \backslash G_0$, and that 
 $u^{G_0}$ fibers over $L_x$ with fiber
$K_0$. Thus $\dim L_x = \dim u^{G_0} - \dim K_0$. On the other hand, since  
\[
{\mathcal O}(v) \cong (G_0 \cap u^{-1}(TAN)u) \backslash G_0,
\]
we know that $u^{G_0}$ fibers over ${\mathcal O}(v)$ with fiber
$ (G_0 \cap u^{-1}(TAN)u) /(G_0 \cap u^{-1}(AN) u)$, which is diffeomorphic to 
$G_0 \cap u^{-1} T u$ by Lemma \ref{lem_diffeo}.   
Thus $\dim u^{G_0} = \dim {\mathcal O}(v)
+ t(v)$, and we have
\[
 \dim L_x = {\rm dim}({\mathcal O}(v))-{\rm dim}(K_0)+t(v).
\]
The formula for the co-dimension of 
$L_x$ in $U/K$ now follows from the facts that the co-dimension of 
${{\mathcal O}}(v)$ 
in $Y$ is
$l(\psi(v)w_bw_0)$ and that $t(v) + \a(v) = \dim T$.

Let $t \in T$. Then $tL_x = L_x$ if and only if there exists $g_0 \in G_0$ 
such that
$tuK_0 = u^{g_0}K_0 \in U/K_0$. By replacing $g_0$ by a product of $g_0$ 
with some $k_0 
\in K_0$, we see that $tL_x = L_x$ if and only if there exists $g_0 \in G_0$
such that $tu = u^{g_0}$, which is equivalent to $bt \in uG_0 u^{-1}$ for some 
$b \in AN$. By Lemma \ref{lem_diffeo}, this is equivalent to 
$t \in T \cap uG_0u^{-1} = T^{\psi(v)\tau}$.
\qed

By \cite[Proposition 1.3.1.3]{Warner}, for every $v \in V$, we can always 
choose $u \in {{\mathcal V}}$
representing $v$ such that $\fg_0 \cap \Ad_{u}^{-1} \fa = 
(\Ad_{u}^{-1} \fa)^\tau \subset \fa^\tau$. When ${{\mathcal O}}(v)$ is 
open in $Y$, 
 $\fg_0 \cap \Ad_{u}^{-1} \fh$ is a
maximally compact Cartan subalgebra of $\fg_0$ \cite{wolf1}, which is unique
up to $G_0$-conjugation. Let $\fh_1$ be any maximally compact
Cartan subalgebra of $\fg_0$ whose vector part $\fa_1$ lies in
$\fa_0 = \fa^\tau$, and let $\fa_{0}^{\prime}$ be any complement of $\fa_1$ in
$\fa_0$. Let $A_{0}^{\prime} = \exp \fa_{0}^{\prime} \subset A_0$.
We have the following corollary of Lemma \ref{lem_leaves-double-cosets}
and Theorem \ref{thm:83}.

\begin{cor}
\label{cor_largest-leaves}
 A symplectic leaf of $\pi_0$ has the largest dimension among all
symplectic leaves if and only if it lies in ${{\mathcal S}}(v)$ 
corresponding to an open $G_0$-orbit ${{\mathcal O}}(v)$.
Such a leaf is diffeomorphic to $ A_{0}^{\prime}N_0$.
\end{cor}

\begin{cor}
\label{cor:when-d-1}
The Poisson structure $\pi_0$
has open symplectic leaves if and only if $\fg_0$ has
a compact Cartan subalgebra. In this case the number of
open symplectic leaves of $\pi_0$ is the same as the 
number of open $G_0$-orbits in $Y$, and each open symplectic
leaf is diffeomorphic to $G_0/K_0$.
\end{cor}

For the rest of this section we assume that $X=U/K_0$ is an irreducible 
Hermitian symmetric space. In this case, there is a parabolic 
subgroup $P$ of $G$ containing $B=TAN$ such that $u_0 K_0 u_{0}^{-1} = U \cap P$
for some $u_0 \in U$. It is proved in \cite{LW} that 
the Poisson structure $\pi_U$ on $U$ projects to a Poisson
structure on $U/(U \cap P)$, which can be regarded as
a Poisson structure on $U/K_0$, denoted by $\pi_\infty$,  
via the $U$-equivariant identification
\[
X=U/K_0 \lrw U/(U\cap P): \, \, 
uK_0 \longmapsto uu_{0}^{-1}(U \cap P). 
\]
Since $(X, \pi_\infty)$ is also $(U, \pi_U)$-homogeneous, the difference
$\pi_0 -\pi_\infty$ is a $U$-invariant bivector field on $X$.
On the other hand, $X$ carried a $U$-invariant symplectic 
structure which is unique up to scalar multiples. Let $\omega_{{\rm inv}}$ be 
such
a symplectic structure, and let $\pi_{{\rm inv}}$ be the corresponding Poisson
bi-vector field. Then since every $U$-invariant bi-vector field 
on $X$ is a scalar multiple of $\pi_{{\rm inv}}$, we have

\begin{lem} There exists $b\in\R$ 
such that $\pi_0=\pi_\infty+b\cdot\pi_{{\rm inv}}$.
\end{lem}

The family of Poisson structures $\pi_\infty+b\cdot\pi_{{\rm inv}}$, $b \in \R$,
has been studied in \cite{KRR}. 
We also remark that when $X$ is Hermitian
symmetric, it is shown in \cite{R-S:orbits} that 
there is a way of parameterizing the $G_0$-orbits in $Y$,
and thus symplectic leaves of $\pi_0$ in $X$,
using only the Weyl group $W$. We refer the interested reader
to \cite[Section 5]{R-S:orbits}. 

\begin{exam}
{\em  
Consider the case $\fg=\sl(2, \C)$, $\fg_0=\sl(2, \R)$. We have  
$U={\rm SU}(2)$, and $K_0$ is the subgroup of $U$ isomorphic to $S^1$ given by:
$$ K_0=\left\{ \left( \begin{array}{cc} \cos t & \sin t \\ -\sin t & \cos t
\end{array} \right), \  \ t\in\R   
\right\} .
$$
The space $X=U/K_0$ can be naturally identified with the Riemann sphere $S^2$ via the map
$$
M=\left( \begin{array}{cc} a & b \\ -{\bar b} & {\bar a} \end{array} \right) \mapsto
z=\frac{-{\rm Im}(a)+\i\cdot{\rm Im}(b)}{{\rm Re}(a)+\i\cdot{\rm Re}(b)},
$$
where $M\in {\rm SU}(2)$ with $|a|^2+|b|^2=1$ and $z$ is a 
holomorphic coordinate on $X\setminus\{ {\rm pt}\}$.   
Then the Poisson structure $\pi_0$ is given by
$$
\pi_0=\i (1-|z|^4)
\frac{\partial}{\partial z}\wedge\frac{\partial}{\partial{\bar{z}}}.
$$
Therefore there are two open symplectic leaves for $\pi_0$, which can be thought 
of as the Northern and the Southern hemispheres. Every point in the Equator, 
corresponding to $|z|=1$, 
is a symplectic leaf as well. It is interesting to notice that the image of a 
symplectic leaf in $U$
given by:  
$$ \frac{1}{\sqrt{1+|z|^2}}\left(  
\begin{array}{cc} z & 1 \\ -1 & {\bar z}    
\end{array} \right) , \ \ z\in\C
$$  
is the union of the Northern and the Southern hemispheres and a point in the Equator.
All three are Poisson submanifolds of $S^2$.   
}
\end{exam}

\begin{rem}
{\em Let ${\mathcal L}$ be the variety of Lagrangian subalgebras of  
$\fg$ with respect to the  
pairing ${\rm Im}\ll \, , \, \gg$, as defined in \cite{EL}.  
Then $G$ acts on  ${\mathcal L}$ by conjugating the subalgebras.  
The variety ${\mathcal L}$  
carries a Poisson structure $\Pi$ defined by the Lagrangian  
splitting $\fg=\fu+(\fa+\fn)$  
such that every $U$-orbit (as well as every $AN$-orbit)  
is a Poisson subvariety of  
$({\mathcal L}, \Pi)$. Consider the point $\fg_0$  
of ${\mathcal L}$ and let $X'$ be the $U$-orbit in  
${\mathcal L}$ through $\fg_0$. Then we have a natural map
$$
{\mathcal J}:\ \   U/K_0 \longrightarrow X'.
$$
The normalizer subgroup of $\fg_0$ in $U$ is not  
necessarily connected but always has $K_0$ as  
its connected component. Thus ${\mathcal J}$ is a finite
covering map. It follows from \cite{EL} that the map ${\mathcal J}$ is 
Poisson with respect to the
Poisson structure $\Pi$ on $X'$.
}
\end{rem}

\section{Invariant Poisson cohomology of $(U/K_0, \pi_0)$.}

Let $\chi^{\fd}(X)$ stand for the graded vector space of the multi-vector  
fields on $X$. Recall that the Poisson coboundary operator, 
introduced by  
Lichnerowicz \cite{Lich}, is given by:
$$
d_{\pi_0}: \ \chi^i(X)\to \chi^{i+1}(X), \ \  
d_{\pi_0}(V)=[\pi_0, V],$$
where $[\cdot, \cdot ]$ is the Schouten bracket of the multi-vector fields
\cite{Kosz}. The Poisson cohomology of $(X, \pi_0)$ is defined to be
the cohomology of the cochain complex $(\chi^{\fd}(X), d_{\pi_0})$ and
is denoted by $H^{\fd}_{\pi_0}(X)$.
By \cite{LuDuke}, the space $(\chi^{\fd}(X))^U$  
of $U$-invariant multi-vector fields on $X$ is closed under $d_{\pi_0}$.
The  cohomology of the cochain sub-complex $((\chi^{\fd}(X))^U, d_{\pi_0})$
is called the $U$-invariant Poisson cohomology of $(X, \pi_0)$ and
we denote it by $H^{\fd}_{\pi_0, U}(X)$. 
We have the following result from \cite[Theorem 7.5]{LuDuke}, adapted to  
our situation $X=U/K_0$, which relates the Poisson cohomology of a Poisson homogeneous  
space with certain relative Lie algebra cohomology. Recall that $G_0$, as
a subgroup of $G$, acts on $U$ by (\ref{eq_G-on-U}), and thus 
$C^{\infty}(U)$ can be viewed as a  
$\fg_0$-module. We also treat $\R$ as the trivial $\fg_0$-module:    

\begin{prop} \cite{LuDuke}
$$
 H^{\fd}_{\pi_0}(X)\simeq H^{\fd}(\fg_0, \fk_0, C^{\infty}(U)), \mbox{  and  }
 H^{\fd}_{\pi_0, U}(X)\simeq H^{\fd}(\fg_0, \fk_0, \R),
$$
\end{prop}

We will compute the cohomology space $H^{\fd}_{\pi_0}(X)$ in a future paper.
The Poisson homology of $\pi_0$ for $X=\C\P^n$ was computed in \cite{kotov}. 
For the $U$-invariant Poisson cohomology, we have

\begin{thm} The $U$-invariant Poisson cohomology of $(U/K_0, \pi_0)$  
is isomorphic to the De Rham cohomology $H^{\fd}(X)$, or, equivalently,  
to the space of $G_0$-invariant differential forms on  
the non-compact dual symmetric space $G_0/K_0$.  
\end{thm}  
\proof By \cite[Corollary II.3.2]{BorWal},  
$H^q(\fg_0, \fk_0, \R)$ is  
isomorphic to $(\wedge^q\fq_0^*)^{\fk_0}$, where $\fq_0$ is the radial part  
in the Cartan decomposition $\fg_0=\fk_0+\fq_0$. This space is isomorphic  
the space of $G_0$-invariant differential $q$-forms on  
the space $G_0/K_0$. Since $\fu=\fk_0+\i\fq_0$,  
and $U$ is compact, we obtain   
$$H^q(\fg_0, \fk_0, \R)\simeq H^q(\fu, \fk_0, \R)\simeq H^q(U/K_0).$$
\qed

\section*{Acknowledgements.}

We thank Sam Evens for many useful discussions.
The first author was partially
supported by NSF grant DMS-0072520. The second author
was partially
supported by NSF(USA) grants DMS-0105195 and DMS-0072551 and by the 
HHY Physical Sciences Fund 
at the University of Hong Kong.
  


\begin{thebibliography}{AAAA}

\bibitem{araki}
Araki, S., On root systems and an infinitesimal classification of
irreducible symmetric spaces, {\em J. Mathematics, Osaka City University},
{\bf 13} (1) (1962), 1-34.

\bibitem{BorWal}{A. Borel and N. Wallach, Continuous cohomology, discrete subgroups,  
and representations of reductive groups. {\it Math. Surveys and Monographs},  
{\bf 67}, A.M.S., 2000.}


\bibitem{EL}{S. Evens and J.-H. Lu, On the variety of Lagrangian subalgebras, I.  
{\it Ann. Scient. \'Ec. Norm. Sup.}, {\bf 34}: 631-668, 2001.}

\bibitem{Fernandes}{R. L. Fernandes, Completely integrable bi-Hamiltonian
systems. {\it Ph.D. Thesis}, U. Minnesota, 1993.}


\bibitem{huckle-wolf}{A. Huckleberry and J. Wolf, Cycle spaces of 
flag domains: a complex geometric viewpoint. arxiv:math.RT/0210445.}

\bibitem{KRR}{S. Khoroshkin, A. Radul, and V. Rubtsov, A family of Poisson structures
on Hermitian symmetric spaces. {\it Comm. Math. Phys.},  
{\bf 152}(2): 299-315, 1993.}

\bibitem{Kosz}{J.-L. Koszul, Crochet de Schouten-Nijenhuis et cohomologie.  
In {\it Math. Heritage of Elie Cartan}, Ast\'erisque, numero hors s\'erie:  
257-271, Soc. Math. France, 1985.}

\bibitem{kotov}{A. Kotov, Poisson homology of $r$-matrix type orbits. I. 
Example of Computation. {\it J. Nonlinear Math. Physics}, {\bf 6}(4): 365-383, 1999.}

\bibitem{Lich}{A. Lichnerowicz, Les vari\'et\'es de Poisson et leurs alg\`ebres  
de Lie associ\'ees. {\it J. Diff. Geom.}, {\bf 12}(2): 253-300, 1977.}

\bibitem{LuDuke}{J.-H. Lu, Poisson homogeneous spaces and Lie algebroids
associated to Poisson actions. {\it Duke Math. J.}, {\bf 86}(2): 261-304, 1997.}

\bibitem{LW}{J.-H. Lu and A. Weinstein, Poisson Lie groups, dressing
transformations, and Bruhat decompositions. {\it J. Diff. Geom.},  
{\bf 31}: 501-526, 1990.}

\bibitem{Matsuki:cartan}{T. Matsuki, The orbits of affine symmetric spaces under
the actionof minimal parabolic subgroups. {\it J. Math. Soc. Japan}, {\bf 31}(2): 
331-357, 1979.}

\bibitem{R-S:orbits}{R. W. Richardson and T. A. Springer,  
Combinatorics and geometry of $K$-orbits on flag manifolds. {\it Contemporary  
Mathematics}, Vol. 153, 109-142, 1993.}


\bibitem{Warner}{G. Warner, Harmonic analysis on semi-simple Lie groups. I. 
{\it Die Gr\"undlehren der mathematischen Wissenschaften}, {\bf 188}, 
Springer-Verlag, 1972.}

\bibitem{wolf1}{J. A. Wolf, The action of a real semisimple Lie group 
on a complex flag manifold, I:
Orbit structure and holomorphic arc components. {\it Bull. Amer. Math. 
Soc.}, {\bf 75} (1969):
1121-1237.}

\end{thebibliography}
\end{document}